\documentclass[11pt]{amsart}
\usepackage{}
\usepackage{amssymb}
\theoremstyle{plain}
\newtheorem{Thm}{Theorem}

\newtheorem{Def}[Thm]{Definition}

\begin{document}

\title[Stability of Cigar soliton]
{Stability of Cigar soliton}

\author{Li Ma, Ingo Witt}
\address{Li Ma, Department of mathematics \\
Henan Normal university \\
Xinxiang, 453007 \\
China} \email{lma@tsinghua.edu.cn}

\address{Ingo Witt, Math. Institut \\
Universitat G\"ottingen \\
Bunsenstr. 3-5, D-37073, G\"ottingen, Germany}

\email{iwitt@uni-math.gwdg.de}

\thanks{The research is partially supported by the National Natural Science
Foundation of China No.11271111 and SRFDP 20090002110019}

\begin{abstract}
In this paper, we use the K\"ahler geometry formulation to
study the global behavior of the Ricci flow on $R^2$. The geometric
feature of our Ricci flow is that it has bounded curvature and finite width. Our aim is to
determine the limiting metric (which corresponds an eternal Ricci
flow) obtained by L.F.Wu. We can use the classification result of
Daskalopoulos-Sesum to give a sufficient condition such that the
limiting metric of L.F.Wu is the metric of cigar soliton.

{ \textbf{Mathematics Subject Classification 2000}: 53C44,32Q20, 58E11}

{ \textbf{Keywords}: Ricci flow, cigar soliton, stability}
\end{abstract}

 \maketitle

\section{Introduction}
In this paper, we study the stability of cigar soliton \cite{H} of
the Ricci flow in $R^2$:
\begin{equation}\label{(1.1)}
\partial_tg(t)=-2Rc(g(t)),
\end{equation}
with the initial metric $g(0)$ of bounded curvature. In the PDE language, (\ref{(1.1)}) is called the Porous-Media equation. The local
existence of the flow on a complete non-compact Riemannian manifold
has been proved by W.Shi in \cite{SH1} (and see \cite{SH2} for the
local derivative estimates via the maximum principle method). Ricci
flow on $R^2$ can be considered as a degenerate diffusion \cite{DK}.
In the works \cite{Hs} and \cite{W}, there are interesting results
about the global Ricci flow on $R^2$. In particular, in \cite{W},
the global Ricci flow on $R^2$ is proved provided the initial
complete metric has bounded curvature and bounded quantity $|d \log
tr_{g_E}g(0)|_{g(0)}$, where $g_E$ is the standard Euclidean metric
on $R^2$ (and we shall always these two conditions are true).
Furthermore, when the initial metric has positive curvature and
finite circumference at infinity, then the flow converges to the
cigar metric at the time infinity. Recall here that the
\emph{circumference at infinity} of the Ricci 2-manifold $(R^2,g)$,
which is defined as
$$
C_\infty(g)=\sup_K\inf_D\{L(\partial D); \text{ for any compact $K$
and open $D$, $K\subset D\subset R^2$}\}.
$$
One may find more references about recent development of Ricci flow
in $R^2$ in \cite{DS} \cite{DS2} \cite{IJ} \cite{JMS}. The aim of
this paper is to extend the above result of L.F. Wu without the
positive curvature assumption.

We now define the width of the evolving metric $g=g(t)$. Let
$F:R^2\to [0,\infty)$ be a proper function F, such that $F^{-1}(a)$
is a compact subset of $R^2$ for every $a\in [0,\infty)$. The width
of F is defined to be the supremum of the lengths of the level
curves of F, namely,
$$w(F) = \sup_c L(F=c).$$
The width $w(g)$ of the metric g is defined to be the infimum
$$
w(g) = \inf_Fw(F).
$$

We now recall some well-known formulae. In dimension two, the Ricci
curvature of the Riemannian manifold $(M^2,g)$ is
$$
Rc(g)=Kg,
$$
where $K=K(g)$ is the Gauss curvature of the metric $g$. Let
$R(g(t))$ be the scalar curvature of the metric $g(t)$. Then
$R(g)=2K(g)$ and the Ricci flow in dimension two is reduced to the
conformal flow
$$
\partial_t g(t)=-R(g(t))g(t).
$$
 We are interested in the Ricci flow on $M^2=R^2$. Given a metric
$g_0$ on $R^2$ and let $g(t)=u(t)g_0$ with $u(t)>0$ being a smooth
function in $R^2$. Then
$$
R(g(t))=u(t)^{-1}(-\Delta_{g_0}\log u(t)+R_0),
$$
where $\Delta_{g_0}$ is the Laplacian operator of the metric $g_0$
with Analytist's sign and $R_0$ is the scalar curvature of $g_0$. We
 simply write $u=u(t)$, $R=R(g(t))$,
etc, when the notations cause no confusion. In this case, the Ricci
flow can be written as

\begin{equation}\label{ricci}
u_t=-Ru, \ \ in
\ \ R^2; \ \ u(0)>0,
\end{equation}
which again can be written as
\begin{equation}\label{ricci+mali}
\partial_tu=\Delta_{g_0}\log u(t)-R_0, \ \ in \ \ R^2, \ \ u(0)>0.
\end{equation}

Taking the time derivative of (\ref{ricci}) we find that
\begin{equation}\label{curv}
R_t=\Delta_g R+R^2.
\end{equation}

Note also that the area element changes at the rate
$$
\frac{d}{dt}dA(g(t))=-RdA(g(t)).
$$

Recall that the global existence of the Ricci flow under suitable
geometric conditions has been studied by L.F.Wu \cite{W}. To handle
the convergence of the flow at time infinity, we need the following
concept of convergence.

\begin{Def}\label{conv}
The Ricci flow $g(t)$ is said to have modified subsequence
convergence, if there exists a 1-parameter family of diffeomorphisms
$\{\Phi(t_j)\}$($t_j\to\infty$) such that there exists a subsequence
(denoted again by $t_j$) such that the sequence $\Phi(t_j)^*g(t_j)$
converges uniformly on every compact set as $t_j\to \infty$.
\end{Def}

We then state the following result of L.F. Wu \cite{W} (see the main
theorem in page 440 in \cite{W} and there is a gap in the proof of the modified subsequence convergence result of in page 447 of L.F.Wu's paper \cite{W}, but we can fix it before we apply it in our case and this is done in section \ref{Preliminary}).

\begin{Thm}\label{wu}
Let $g(t) = e^{\check{u}(t)}g_E$ be a solution to (\ref{(1.1)}) such
that $g(0)= e^{\check{u}_0}g_E$ is a complete metric with bounded
curvature and $|\nabla \check{u}_0|$ is uniformly bounded on $R^2$.
Then the Ricci flow has modified subsequence convergence as $t_j\to
\infty$ with the limiting metric $g_\infty$ being complete metric on
$R^2$. Furthermore, the limiting metric is the cigar soliton if
$C_\infty(g(0))<\infty$ and $g(0)$ has positive curvature.
\end{Thm}

We point out that the diffeomorphisms $\Phi(t_j)$ used in Theorem
\ref{wu} are of the special form
$$
\Phi(t)(a,b)=(e^{\frac{-\check{u}(x_0,t)}{2}}a,e^{\frac{-\check{u}(x_0,t)}{2}}b)=(x_1,x_2)=x,
$$
where $x_0=(0,0)$. The important fact for these diffeomorphisms is
that $$ |\nabla_{g(t)} f(x, t)|=|\nabla_{\Phi(t)^*g(t)} f((a, b),
t)|.
$$ for any smooth function $f$ and $x=\Phi(t)(a,b)$. Note that in
the convergent part of Theorem \ref{wu}, the limiting metric
corresponds an eternal Ricci flow $g_\infty(t)$ on $R^2$. In fact,
it is proved in Theorem 2.4 in \cite{W} that the curvature and
$|\nabla \check{u}|^2$ of the global flow are uniformly bounded.
According to the classification result \cite{DS} \cite{DS2}, the
limiting flow is cigar soliton provided the width of $g_\infty$ is
finite, i.e., $w(g_\infty)<\infty$. Our aim is to show this is true
in certain circumstances.

We shall assume that there is a potential function $f(0)$ for the
initial metric $g(0)$ in the sense that $R(0)=\Delta_{g(0)}f(0)$ in
$R^2$ (see \cite{BM} for related results). It is clear that the potential function plus any constant is
still a potential function. For the cigar metric
$g_c=\frac{dx^2+dy^2}{1+x^2+y^2}$, we have its potential function
$f_0=\log (1+x^2+y^2)$ in $R^2$. The cigar metric has finite width,
i.e., $w(g_c)<\infty$. We shall take $g_0$ the cigar metric in our
new result.

We shall prove the following result.

\begin{Thm}\label{cigar} Assume that the Riemannian metric
$g(0):=u(0)g_c$ on $R^2$ has bounded curvature, where $g_c$ is the metric of cigar soliton on $R^2$ and $u(0)$ is positive
smooth function in $R^2$ such that both $|\log u(0)|+|d \log u(0)|_{g(0)}$ is uniformly bounded on $R^2$
(so $g(0)$ has finite width). Suppose that the potential function
$f(0)$ of $g(0)$ satisfies that $f_0-f(0)\in L^\infty(R^2,g(0))$.
Then the Ricci flow with the initial metric $g(0)$ has global
solution with its limiting metric at $t=\infty$ the cigar soliton
metric.
\end{Thm}

Our result is different from the result of Hsu \cite{Hs}, which is also
recoded in \cite{IJ} (see also \cite{M}). The result above can be considered as the
stability property of cigar soliton. The global existence of the
Ricci flow (and the limiting metric) follows from the work of L.Wu
\cite{W}. We just need to show that the width $w(g(t))$ is uniformly
bounded. We shall use the K\"ahler-Ricci flow formulation \cite{JMS}
and the maximum principle trick \cite{SH1} to prove theorem above.
One of the important observation in our argument of the result above
is that for the function $v=f(0)-f$, it satisfies that
\begin{equation}\label{lambda}
(\partial_t-\Delta)v=-\Delta f(0)=-u^{-1}u(0)R(0), \ in \ R^2\times
(0,\infty)
\end{equation}
with the initial data $ v(0)=0$ in $ R^2$, where $f$ is the Ricci potential of $g$ such that $R=\Delta_gf$ (see (\ref{f-def}) in next section).

 The structure of this paper is as follows. In the section
\ref{Preliminary}, we recall K\"ahler geometry formulation of the
2-d ricci flow and prove some a priori estimates. Then we fix a gap in the proof of the modified subsequence convergence result of in page 447 of L.F.Wu's paper \cite{W}.  In the section
\ref{proof of theorem}, we give the proof of theorem \ref{cigar}.
Some well known facts about cigar metric is recalled in the last
section.

\section{K\"ahler-Ricci flow formulation and a priori bounds}\label{Preliminary}

We shall use the K\"ahler-ricci flow formulation to study the Ricci
flow on $R^2$. We shall consider the Ricci flow (\ref{ricci}) as the
K\"ahler-Ricci flow by setting
$$
g_{i\bar{j}}={g_0}_{i\bar{j}}+\partial_i\partial_{\bar{j}}\phi,
$$
where $\phi=\phi(t)$ is the K\"ahler potential of the metric $g(t)$
relative to the metric $g_0$. Note that $$
g(0)_{i\bar{j}}={g_0}_{i\bar{j}}+\partial_i\partial_{\bar{j}}\phi_0,
$$
In this situation, the Ricci flow can be written as
\begin{equation}\label{kahler}
\partial_t\phi=4\log \frac{{g_0}_{1\bar{1}}+\phi_{1\bar{1}}}{{g_0}_{1\bar{1}}}-4f_0, \ \ \phi(0)=\phi_0,
\end{equation}
where $f_0$ is the potential function of the metric $g_0$ in the
sense that $R(g_0)=\Delta_{g_0}{f_0}$ in $R^2$. Such a potential
function has been introduced by R.Hamilton in \cite{H}. We remark
that the initial data for the evolution equation (\ref{kahler}) is
$\phi(0)$ which is non-trivial. Note also that we have used the
factor $4$ in right side of (\ref{kahler}) which can allow us to use
the usual Laplacian operator Riemannian geometry, for otherwise, we
need to use the normalized Laplacian in K\"ahler geometry.

Let
\begin{equation}\label{f-def}
f=-\partial_t\phi.
\end{equation}
Then, taking the time derivative of (\ref{kahler}), we have
\begin{equation}\label{potential}
\partial_tf=\Delta_g f, \ \ f(0)=-\partial_t\phi(0).
\end{equation}
By this we can easily get (\ref{lambda}). The important fact for us
is that
\begin{equation}\label{scalar}
\Delta_g f=R.
\end{equation}

To prove this, we recall some well-known results. Recall that any
2-d Riemannian manifold is a 1-d K\"ahler manifold. On a K\"ahler
manifold of complex dimension $n$, the metric in coordinate
expression is
$$
g=g_{i\bar{j}}dz^idz^{\bar{j}}
$$
and the Ricci form is
$$
\rho=\frac{\sqrt{-1}}{2}R_{i\bar{j}}dz^i\wedge
dz^{\bar{j}}=-\sqrt{-1}\partial_i\partial_{\bar{j}}Gdz^i\wedge
dz^{\bar{j}}
$$
with
$$
G=:G(g)=\log det (g_{i\bar{j}}).
$$
Then the Ricci curvature of $g$ is $$
R_{i\bar{j}}=-2\partial_i\partial_{\bar{j}}G
$$
and its scalar curvature is
$$
R(g)=2g^{i\bar{j}}R_{i\bar{j}}=-\Delta_g G,
$$
where $(g^{i\bar{j}})$ is the inverse of the matrix $(g_{i\bar{j}})$
and $\Delta_g=4g^{i\bar{j}}\partial_i\partial_{\bar{j}}$ is the
Laplacian operator of the metric $g$. In the complex dimension one
case, we have
$$
Kg_{1\bar{1}}=R_{1\bar{1}}=-2\partial_1\partial_{\bar{1}}G,
$$
where $K$ is the Gauss curvature of the metric $g$. Then the Ricci
potential function for the metric $g$ is $f=-G$.

Here we compute an example for the help of understanding. For the cigar metric $$ g=g_c:=\frac{|dz|^2}{1+|z|^2},
$$
on $R^2$, we have
$$
g_{1\bar{1}}=\frac{1}{1+|z|^2}, \ \ g^{1\bar{1}}={1+|z|^2}
$$
and
$$
G=-\log (1+|z|^2).
$$
Then $R=\frac{4}{1+|z|^2}$, $K=\frac{2}{1+|z|^2}$, and
$$
Kg_{1\bar{1}}=\frac{2}{(1+|z|^2)^2}=2\partial_1\partial_{\bar{1}}\log
(1+|z|^2).
$$
Hence
$$
f(z)=-G=\log (1+|z|^2)=f_c(z)
$$
is the Ricci potential of the metric $g$.

 We now prove (\ref{scalar}). Note that
$$
Rg_{1\bar{1}}=2R_{1\bar{1}}=-4\partial_1\partial_{\bar{1}}G.
$$
Applying $\partial_1\partial_{\bar{1}}$ to (\ref{kahler}) we obtain
that
$$
-4\partial_1\partial_{\bar{1}}f=4\partial_1\partial_{\bar{1}}G-4\partial_1\partial_{\bar{1}}G(g_0)-4\partial_1\partial_{\bar{1}}f_0.
$$
Taking the trace with respect to $g_{1\bar{1}}$ we have
$$
-\Delta_g f=-R+\frac{1}{u}(R_0-\Delta_{g_0} f_0)=-R,
$$
which is (\ref{scalar}). Here we have used
 that
\begin{equation}\label{fact1}
\Delta_g:=\Delta =u^{-1}\Delta_{g_0}
\end{equation}
and
\begin{equation}\label{fact2}
R_0=\Delta_{g_0}f_0.
\end{equation}

Using (\ref{fact1}) and (\ref{fact2}), we can write
(\ref{ricci+mali}) as
\begin{equation}\label{key}
(\log u-f_0)_t=\Delta (\log u-f_0),  \ \ in \ \ R^2.
\end{equation}

We define $w$ by
$w=\log u +f-f_0$.
Then we have
 $$ u=e^{w-f+f_0}=e^{w+v+f_0-f(0)}.
$$
 Then adding (\ref{potential}) and (\ref{key}), we have
\begin{equation}\label{max1}
\partial_tw=\Delta_g w, \ \ in  \ \  R^2
\end{equation}
with the initial data
$$
w(0)=\log u(0)-f_0+f(0), \ \ in  \ \  R^2,
$$
which is a bounded function in $R^2$ by our assumption.

Note that by (\ref{ricci}), (\ref{potential}) and (\ref{scalar}), we have
$$
w_t=(\log u)_t+f_t=-R+\Delta_g f=0,  \ \ in \ \ R^2.
$$
Then $w(x,t)=w(x,0)$ for any $x\in R^2$ and $t>0$.

We now can give an outline of the proof of Theorem \ref{cigar}: We
shall show in next section that $v=f(0)-f$ is uniformly upper
bounded. Once this is done, we then can conclude that
\begin{equation}\label{ma}
u(t)=e^{w+v-f(0)+f_0}
\end{equation}
is uniformly upper bounded, which makes the metric $g(t)$ be
uniformly upper bounded by a scale of the cigar soliton and the
limiting metric $g_\infty$ have finite width $w(g_\infty)$, which is
the cigar soliton, by using the classification result of ancient
solutions due to Daskalopoulos and Sesum \cite{DS} \cite{DS2}. This
will completes the proof of Theorem \ref{cigar}.

We now fix the gap in the proof of the modified subsequence convergence result of in page 447 of L.F.Wu's paper \cite{W} used in our case. Recall that
$$
g=g(t)=u(t)g_0 = e^{\check{u}(t)}g_E
$$
where $\check{u}(0)=\log u_0-f_0$. Note that by (2.2) in \cite{W}, $\check{u}=\check{u}(t)$ satisfies on $R^2$ that
$$
\partial_t\check{u}=\Delta_{g}\check{u}=-R.
$$
By the uniform upper bound assumption about $u(t)$,  $\check{u}(t)$ is uniformly upper bounded in $R^2$.
Actually,  the uniform upper bound for $\check{u}(t)$ is given by
$$
\sup_{R^2} \check{u}(t)\leq \sup_{R^2} \check{u}(0)=:C
$$
for all time. 

By our assumption, $g(0)$ is uniformly equivalent to the cigar metric and $|\nabla\check{u}(0)|_{g(0)}$ is uniformly bounded. By Theorem 2.4 in \cite{W} we know that there are uniform constant $k_0$ and $D_0$ such that
$$
|R|\leq k_0+D_0, \ \ \ \ \text{and} \ \ \ |\nabla \check{u}(t)|_{g(t)}^2\leq D_0
$$
for all time. The latter bound can be written as
\begin{equation}\label{star}
|d \check{u}(t)|_{g_E}^2\leq D_0 e^{\check{u}(t)}\leq D_0 e^C
\end{equation}
for all time. Applying the diffeomorphism $\Phi(t)$ as in page 447 in \cite{W} we may assume that $\check{u}(0,t)=0$ and denoted by $\check{g}(t)=\Phi(t)^*g(t)$. By (\ref{star}), uniform estimates of the higher derivatives of $\check{u}(t)$ and curvatures of the metric $g(t)$ in Corollary 2.5 in \cite{W} and by Arzela-Ascoli theorem, for any sequence of times going to infinity there exists a subsequence $\{t_j\}$ such that $\check{u}(t_j)$ converges uniformly to a Lipschtz function $\check{u}_\infty$ on every compact subset ( and correspondingly the metric sequence $\check{g}(t_j)$ converges uniformly to a Lipschtz function $\check{g}_\infty=e^{\check{u}_\infty}g_E$ on every compact subset). This fills in the gap in L.F.Wu's result in our case.

\section{main estimates and the proof of Theorem \ref{cigar}} \label{proof of theorem}

In this section we shall give the uniform upper bound of $v$, which implies the uniform upper bound of $u$ and it
is enough for us to conclude the desired geometric width bound of the limiting metric $g_\infty$ at time infinity of the Ricci flow.

Note that in the statement of Theorem \ref{cigar}, we have $g_0=g_c$ and $g(0)=u_0g_c$, where $g_c$ is the
metric of the cigar soliton and $u_0>0$ is a bounded positive
function in $R^2$. For convenient of the readers, we recall some
facts about the cigar soliton in $R^2$ in the appendix. Good
references for this metric are \cite{H} and \cite{BK}. The crucial
fact for the cigar metric is that $$ C_\infty(g_c)<\infty.
$$

Note that $$ R(0)=u_0^{-1}(-\Delta_{g_c}\log u_0+R_c).
$$
Then for $u=u(t)$,
$$
R(0)=u_0^{-1}(-u\Delta_{g}\log u_0+R_c).
$$
Inserting this into (\ref{lambda}) we get that
$$(\partial_t-\Delta)v=\Delta_{g}\log u_0-u^{-1}R_c.
$$
Recall that $v=v(t)=f(0)-f$ and $v(0)=0$. By (\ref{scalar}) and (\ref{potential}) we have
$$
-v_t=(f-f(0))_t=f_t=R.
$$
Then,
$$
v:=v(t)=-\int_0^t R,
$$
which is a uniformly bounded function on any finite interval $[0,T]$. 

Let $h=v+\log u_0$. Then we have
$$
(\partial_t-\Delta)h=-u^{-1}R_c<0, \ in \ R^2\times (0,\infty),
\ h(0)=\log u_0.
$$
By assumption that $\log u_0$ is uniformly bounded on $R^2$, so is $h$ on $R^2$ for any finite interval $[0,T]$. By the maximum principle (Theorem 1.1 in \cite{W}), we get the uniform bound for
$h$ on $R^2$ for all time. Then we know that $v=f(0)-f=h-\log u_0$ is uniformly bounded in
$R^2$ for all time, and hence $g(t)$ is uniformly controlled by the cigar metric
$g_c$ up to a factor. Therefore, the width
$$
w(g(t))<\infty
$$
is uniformly bounded, and then
$$
w(g_\infty)<\infty.
$$
Then, the limiting metric is a cigar metric and we have completed
the proof of Theorem \ref{cigar}.

\section{Appendix: cigar soliton}
Since the stability of cigar metric is our main topic here, we
prefer to review a little more about it.
 Hamilton's cigar soliton is a special solution to the
Ricci flow equation of metrics $g(t)$ in $R^2$:
$$ g_t=-Rg, \ \ in \ \ R^2.
$$
Namely, it is the one-parameter family of complete Riemannian
metrics of the form
$$
g_c(t)=\frac{dx^2+dy^2}{e^{4t}+x^2+y^2}
$$
in $R^2$. Define the diffeomorphism
$$\phi_t(a,b)=(e^{2t}a,e^{2t}b)=(x,y)
$$
and the metric
$$
g_c(a,b):=\phi_t^{*}g_c(t)=\frac{da^2+db^2}{1+a^2+b^2},
$$
which is the so called the cigar metric.

We now re-write the cigar metric in polar coordinates $(r,\theta)$
with $r=\sqrt{a^2+b^2}$ such that
$$
g_c=\frac{dr^2+r^2d\theta^2}{1+r^2}.
$$
One can compute (see pages 24-28 in \cite{BK}) that the scalar curvature of $g_c$ is
$$
R_c=\frac{4}{1+r^2}
$$
and its area element is
$$
dv_c=\frac{rdrd\theta}{1+r^2}.
$$

Let $z=a+\sqrt{-1}b$. Then we have
$$
g_c=\frac{|dz|^2}{1+|z|^2},
$$
which is the K\"{a}hler metric on $R^2:=\mathcal{C}$. Define the
parameter $s$ such that
$$
ds=\frac{dr}{\sqrt{1+r^2}}.
$$
Then we have
$$
s=\log (r+\sqrt{1+r^2})=arcsinh r
$$
and
$$
g_c=ds^2+tanh^2sd\theta^2.
$$
Then
$$
R_c=4(cosh s)^{-2}.
$$

We now try other method of finding the Ricci potential $f_0$ of the
metric $g_c$, i.e.,
$$\frac{R_c}{2}g_c=\nabla ^2 f_0.
$$
We shall look for $f$ being the radial function $f_0=f_0(s)$. Define
the local normal frame $e_1=\frac{\partial}{\partial s}$ and $
e_2=\frac{1}{tanh s}\frac{\partial}{\partial \theta}.
$ Then
$$
\frac{R_c}{2}=\frac{R_c}{2}g_c(e_1,e_1)=\nabla^2f(e_1,e_1)=f_0^{''}(s)
$$
and
$
\frac{R_c}{2}=\frac{R_c}{2}g_c(e_2,e_2)=-(\nabla_{e_2}e_2) f_0.
$

By these relations we find that
$
f_0=2\log cosh s.
$
Let $g_c=w_0(dx^2+dy^2)$ in $R^2$, where
$$
w_0(x,y)=\frac{1}{1+x^2+y^2}=(cosh s)^{-2}.
$$
The key observation is that we have the relation
$$
\log w_0+f_0=0.
$$
This fact plays a role in our formulation of Theorem \ref{cigar}.

{\bf Acknowledgement}. The author thanks the unknown referees very much for useful suggestions, in particular for pointing out the gap in the proof of the modified subsequence convergence result of in page 447 of L.F.Wu's paper \cite{W}.

\end{document}